\documentclass[a4paper,reqno,12pt]{amsart}
\usepackage{enumerate}
\usepackage{amssymb}
\usepackage{ifthen}
\usepackage{graphicx}
\usepackage{color,xcolor}
\usepackage{mathrsfs}
\numberwithin{equation}{section}
\newtheorem{thm}{Theorem}
\newtheorem{theorem}{Theorem}
\newtheorem{lem}{Lemma}
\newtheorem{cor}{Corollary}
\newtheorem{defn}{Definition}
\newtheorem{example}{Example}

\newtheorem{claim}{Claim}

\newtheorem{conjecture}{Conjecture}
\newtheorem{conj}{Conjecture}

\newtheorem{prob}{Problem}

\newenvironment{rem}{%
\bigskip
\noindent \textsl{{\sl Remark. }}}{\bigskip}
\newenvironment{rems}{%
\bigskip
\noindent \textsl{{\sl Remarks. }}}{\bigskip}

\newenvironment{pf}[1][]{%
 \vskip 1mm
 \noindent
 \ifthenelse{\equal{#1}{}}%
  {{\slshape Proof. }}%
  {{\slshape #1.} }%
 }%
{\qed\medskip}

\newcounter{alphabet}
\newcounter{tmp}

\makeatletter
\newcommand{\Ref}[1]{\@ifundefined{r@#1}{}{\setcounter{tmp}{\ref{#1}}\Alph{tmp}}}
\makeatother

\DeclareMathOperator{\RE}{Re}

\def\be{\begin{equation}}
\def\ee{\end{equation}}
\def\bes{\begin{equation*}}
\def\ees{\end{equation*}}

\newcommand{\bee}{\begin{enumerate}}
\newcommand{\eee}{\end{enumerate}}

\newcommand{\blem}{\begin{lem}}
\newcommand{\elem}{\end{lem}}
\newcommand{\bthm}{\begin{thm}}
\newcommand{\ethm}{\end{thm}}
\newcommand{\bcor}{\begin{cor}}
\newcommand{\ecor}{\end{cor}}
\newcommand{\beg}{\begin{example}}
\newcommand{\eeg}{\end{example}}
\newcommand{\begs}{\begin{examples}}
\newcommand{\eegs}{\end{examples}}
\newcommand{\bdefe}{\begin{defn}}
\newcommand{\edefe}{\end{defn}}
\newcommand{\bprob}{\begin{prob}}
\newcommand{\eprob}{\end{prob}}
\newcommand{\bques}{\begin{ques}}
\newcommand{\eques}{\end{ques}}
\newcommand{\bei}{\begin{itemize}}
\newcommand{\eei}{\end{itemize}}

\newcommand{\bde}{\begin{deter}}
\newcommand{\ede}{\end{deter}}
\newcommand{\bca}{\begin{case}}
\newcommand{\eca}{\end{case}}
\newcommand{\bcl}{\begin{claim}}
\newcommand{\ecl}{\end{claim}}
\newcommand{\bcon}{\begin{conj}}
\newcommand{\econ}{\end{conj}}
\newcommand{\bcons}{\begin{conjs}}
\newcommand{\econs}{\end{conjs}}
\newcommand{\bprop}{\begin{propo}}
\newcommand{\eprop}{\end{propo}}
\newcommand{\br}{\begin{rem}}
\newcommand{\er}{\end{rem}}
\newcommand{\brs}{\begin{rems}}
\newcommand{\ers}{\end{rems}}
\newcommand{\bo}{\begin{obser}}
\newcommand{\eo}{\end{obser}}
\newcommand{\bos}{\begin{obsers}}
\newcommand{\eos}{\end{obsers}}
\newcommand{\bpf}{\begin{pf}}
\newcommand{\epf}{\end{pf}}
\newcommand{\ba}{\begin{array}}
\newcommand{\ea}{\end{array}}
\newcommand{\beq}{\begin{eqnarray}}
\newcommand{\beqq}{\begin{eqnarray*}}
\newcommand{\eeq}{\end{eqnarray}}
\newcommand{\eeqq}{\end{eqnarray*}}

\begin{document}
\bibliographystyle{amsplain}
\title[Harmonic pre-Schwarzian and its applications]
{Harmonic pre-Schwarzian and its applications}


\author[G. Liu]{Gang Liu}
\address{G. Liu, College of Mathematics and Statistics
 (Hunan Provincial Key Laboratory of Intelligent Information Processing and Application),
Hengyang Normal University, Hengyang,  Hunan 421002, China.}
\email{liugangmath@sina.cn}

\author[S. Ponnusamy]{Saminathan Ponnusamy
}
\address{S. Ponnusamy, Stat-Math Unit,
Indian Statistical Institute (ISI), Chennai Centre,
110, Nelson Manickam Road,
Aminjikarai, Chennai, 600 029, India.}
\email{samy@isichennai.res.in, samy@iitm.ac.in}



\begin{abstract}
The primary aim of this article is to extend certain inequalities concerning the pre-Schwarzian derivatives from
the case of analytic univalent functions to that of univalent harmonic mappings defined on certain domains.
This is done in two different ways. One of the ways is to connect with a conjecture on the univalent harmonic mappings.
Also, we improve certain known results on the majorization of the Jacobian of functions in the
affine and linear invariant family of sense-preserving harmonic mappings. This is achieved as an application of
a corresponding distortion theorem in terms of the harmonic pre-Schwarzian derivative.
\end{abstract}

\maketitle \pagestyle{myheadings}
\markboth{G. Liu and S. Ponnusamy}{Harmonic pre-Schwarzian and its applications}

\section{Introduction}  \label{sec1}
Throughout the paper, let $\mathbb {D}=\{z:\,|z|<1\}$, $\mathbb{H}=\{z:\,\RE z>0\}$ and $\Delta=\overline{\mathbb{C}}\backslash\overline{\mathbb{D}}$
denoted the unit disk, the right half-plane and the exterior of the closed unit disk, respectively.
Suppose $f$ (resp. $g$) is  analytic in $\mathbb{D}$ (resp. $\mathbb{H}$) and
$F$ is analytic in $\Delta\backslash\{\infty\}$ with a simple pole at $z=\infty$.
If $f$, $g$ and $F$ are univalent in  $\mathbb{D}$, $\mathbb{H}$ and $\Delta$, respectively,
then we have
\be
\sup_{z\in\mathbb{D}}(1-|z|^2)\left|f''(z)/f'(z)\right|\leq6, \label{equ1.1}
\ee
\be
\sup_{z\in\mathbb{H}}2\RE z\left|g''(z)/g'(z)\right|\leq6, \label{equ1.2}
\ee
and
\be
\sup_{z\in\Delta}(|z|^3-|z|)\left|F''(z)/F'(z)\right|\leq6. \label{equ1.3}
\ee
The constant 6 is sharp in all the three cases.
Inequalities \eqref{equ1.1} and \eqref{equ1.2} are obtained as a consequence of Bieberbach's distortion theorem.
However, the estimate for $F$ is deeper and is established in \cite{avh} as a consequence of Goluzin's inequality \cite[p. 139]{gol}.

Osgood \cite{osg} generalized \eqref{equ1.1} and \eqref{equ1.2} to an arbitrary simply connected domain by means of {\it hyperbolic metric}.
In what follows, $D\subset\mathbb{C}$ is a domain with at least two boundary points.
By the uniformization theorem, the hyperbolic metric  $\lambda_{D}(z)$  is induced by
\begin{equation}   \label{equ1.4}
\lambda_{D}(f(z))|f'(z)|=\lambda_{\mathbb{D}}(z)=1/(1-|z|^2),\quad z\in \mathbb{D},
\end{equation}
where $f:\,\mathbb{D}\rightarrow D$ is a (universal) covering mapping onto $D$.
Especially, if $D$ is simply connected, then $f$ is a conformal mapping of $\mathbb{D}$ onto $D$.
The definition is independent of the choice of the covering mapping of $\mathbb{D}$ onto $D$
since two covering mappings of $D$ differ only by a conformal self-mapping of $\mathbb{D}$.
Its Gaussian curvature is $-4$.

\begin{thm} { \rm(\cite[Theorem~1]{osg})} \label{thmA}
Let $f$ be a univalent analytic function in a simply connected domain $D\subset\mathbb{C}$.
Then we have the sharp inequality
\begin{equation}\label{equ1.5}
\lambda_{D}^{-1}(z)|f''(z)/f'(z)|\leq8,\quad z\in D.
\end{equation}
\end{thm}

Furthermore, Osgood \cite{osg} stated that \eqref{equ1.5} does not hold for  arbitrary domains.
For instance, let $D=\mathbb{D}^*:=\mathbb{D}\backslash\{0\}$ and $f(z)=1/z$.
The basic theory of hyperbolic metric and computation shows that
$$\lambda_{\mathbb{D}^*}(z)=1/(2|z|\log(1/|z|))\quad \mbox{and}\quad
\sup_{z\in\mathbb{D}^*}\lambda_{\mathbb{D}^*}^{-1}(z)|f''(z)/f'(z)|=\infty.
$$
Based on this, he investigated  sufficient and necessary conditions to demonstrate  that such type of inequality
does hold in a multiply connected domain. 

\begin{thm} {\rm(\cite[Theorem~2]{osg})} \label{thmB}
Let $D\subset\mathbb{C}$ have at least two boundary points.
There exists a constant $a$ such that $\lambda_{D}^{-1}(z)|f''(z)/f'(z)|\leq a$ in $D$
for all univalent analytic functions in $D$ if and only if there exists a positive constant $c$ such that
\begin{equation}\label{equ1.6}
\lambda_D(z)d(z,\partial D)\geq c,\quad z\in D.
\end{equation}
Here $d(z,\partial D)$ denotes the Euclidean distance from $z$ to the boundary $\partial D$ of $D$.
\end{thm}

For the characterization of domains satisfying \eqref{equ1.6}, see \cite[Section~5]{osg}.
The operator $f\mapsto P_f:=f''/f'$, is called the {\it pre-Schwarzian derivative} of $f$ when $f$
is locally univalent and analytic in $D$.
The {\it pre-Schwarzian norm} of $f$ in $D$ is defined as $||P_f||_D=\sup_{z\in D}\lambda^{-1}_D(z)|P_f(z)|$.
Recently, Hern\'{a}ndez and Mart\'{\i}n  \cite{HM2015} extended these  notions to  any locally univalent harmonic mapping
in a domain (see also \cite{CDO}).
In view of this, in Section \ref{sec3}, we  generalize those inequalities
to complex-valued  univalent harmonic mappings.
Our results are based on a distortion theorem concerning the {\it harmonic pre-Schwarzian derivative}
about {\it linear and affine invariant family} of harmonic mappings.
However, the sharp results should be resorted to Clunie and Sheil-Small's conjecture (\cite{CS}) in terms {\it specified order}
of the family $\mathcal{S}_H$ of normalized univalent harmonic mappings in the unit disk. In Section \ref{sec4}, we build some sharp inequalities on a class of analytic pre-Schwarzian derivatives
and the harmonic pre-Schwarzian derivatives,  and propose a conjecture
about the relationship between univalent analytic functions and univalent harmonic mappings.
We use these combinations to obtain similar results, as that of Section \ref{sec3}, when the domain is simply connected.
Furthermore, the proposed conjecture implies the former conjecture of Clunie and Sheil-Small.
In Section \ref{sec5}, we improve the  corresponding results of \cite[Section~3]{sch} by applying our distortion theorem
and the results on {\it majorization-subordination} theory of  {\it universal linear invariant family} of analytic functions.

\section{Background and Preliminaries} \label{sec2}

\subsection{Univalent harmonic mappings}  \label{sec2.1}

A complex-valued function $f$ is called a {\it harmonic mapping}
if it satisfies the  Laplace equation $\Delta f=4f_{z\overline{z}}=0$.
In a simply connected domain $D$, every  harmonic mapping $f$
has a {\it decomposition} $f=h+\overline{g}$, where $h$ and $g$ are analytic functions in $D$.
However, in a multiply connected domain, the representation $f=h+\overline{g}$ is valid
locally but may not have a single-valued global extension.

According to Lewy's theorem \cite{lew},  a  harmonic mapping $f$ of the form $f=h+\overline{g}$ is locally univalent in  a domain $D$
if it is sense-preserving, i.e., its Jacobian $J_f=|h'|^2-|g'|^2$ is positive in  $D$
so that  its dilatation $\omega_f$ defined by $\omega_f=g'/h'$ has the property that $|\omega_f(z)|<1$ in $D$.
In the study of univalent harmonic mappings, it is convenient to consider the class $\mathcal{S}_H$
of all sense-preserving and univalent harmonic mappings $f=h+\overline{g}$ in the unit disk $\mathbb{D}$ with the normalizations $h(0)=h'(0)-1=g(0)=0$.
The class $\mathcal{S}_H$ is not compact whereas $\mathcal{S}_H^0:=\{f\in\mathcal{S}_H:\,g'(0)=0\}$ is compact
(see \cite[p. 78]{dur2004}).

Let $\mathcal{K}_H$ (resp. $\mathcal{C}_H$) be the set of all  convex
 (resp. close-to-convex) harmonic mappings from $\mathcal{S}_H$.
Let $\mathcal{K}_H^0:=\mathcal{K}_H\cap\mathcal{S}_H^0$ and $\mathcal{C}_H^0:=\mathcal{C}_H\cap\mathcal{S}_H^0$.
Clunie and Sheil-Small \cite{CS} constructed
the {\it  harmonic half-plane mapping $L$} and the {\it harmonic Koebe function $K$} defined by
$$L(z)=\frac{2z-z^2}{2(1-z)^2}
+\overline{\frac{-z^2}{2(1-z)^2}}\quad \mbox{and}
\quad K(z)=\frac{z-\frac{1}{2}z^2+\frac{1}{6}z^3}{(1-z)^3}
+\overline{\frac{\frac{1}{2}z^2+\frac{1}{6}z^3}{(1-z)^3}},
$$
respectively. These functions  play the role extremal in many problems of harmonic mappings.
Note that $L\in\mathcal{K}_H^0$ and $K\in\mathcal{C}_H^0$.
Basic information about harmonic mappings may be obtained from the monograph of Duren \cite{dur2004} and the recent survey \cite{PR}.

\subsection{Affine and linear invariant families}  \label{sec2.2}
Let $\mathcal{F}$  be a family of sense-preserving harmonic mappings $f=h+\overline{g}$ in $\mathbb{D}$,
normalized by $h(0)=g(0)=h'(0)-1=0$.
$\mathcal{F}$  is said to be a {\it linear invariant family} (LIF) if for each $f\in\mathcal{F}$,
\bes
K_\varphi(f(z))=\frac{f(\varphi(z))-f(\varphi(0))}{\varphi'(0)h'(\varphi(0))}\in\mathcal{F}
\quad \forall ~\varphi\in{\rm Aut}(\mathbb{D}),
\ees
and $\mathcal{F}$ is called an {\it affine invariant family} (AIF) if for each $f\in\mathcal{F}$,
\bes
A_\varepsilon(f(z))=\frac{f(z)+\varepsilon\overline{f(z)}}{1+\varepsilon g'(0)}\in\mathcal{F}
\quad \forall ~\varepsilon\in\mathbb{D}.
\ees
Here $K_\varphi(f)$ and $A_\varepsilon(f)$ are called Koebe and affine transforms of $f$, respectively.
We say that $\mathcal{F}$ is an affine and linear invariant family (ALIF)
if it is both LIF and AIF.
For example, each of $\mathcal{S}_H$, $\mathcal{K}_H$ and $\mathcal{C}_H$ is an ALIF.
The {\it order} of ALIF $\mathcal{F}$, defined by
$$\alpha(\mathcal{F})=\sup_{f\in\mathcal{F}}|a_2(f)|=\frac{1}{2}\sup_{f\in\mathcal{F}}|h''(0)|,
$$
plays an important role in the study of harmonic mappings, since the appearance of the pioneering work of Clunie and Sheil-Small (see \cite{CS}).
In 2007, the notion of {\it specified order} of ALIF $\mathcal{F}$ was introduced in \cite{gra2007} as
$$\alpha_0(\mathcal{F})=\sup_{f\in\mathcal{F}^0}|a_2(f)|=\frac{1}{2}\sup_{f\in\mathcal{F}^0}|h''(0)|,
$$
where $\mathcal{F}^0=\{f=h+\overline{g}\in\mathcal{F}: \,g'(0)=0\}$.
Note that $1/2\leq\alpha_0(\mathcal{F})\leq\alpha(\mathcal{F})\leq\alpha_0(\mathcal{F})+1/2$
and $\alpha(\mathcal{F})\geq1$ for any ALIF $\mathcal{F}$ (see  \cite{gra2016}).
Also the specified order $\alpha_0(\mathcal{F})$ of a given ALIF $\mathcal{F} $ coincides with the {\it new order}
of LIF defined in \cite{SSS}.  Please refer to \cite{gra2012} for further details about $\alpha_0(\mathcal{F})$.

It follows from \cite{CS,WLZ} that
$$\alpha_0(\mathcal{K}_H)=3/2,\quad \alpha_0(\mathcal{C}_H)=5/2,\quad \alpha(\mathcal{K}_H)=2 \quad \text{and} \quad \alpha(\mathcal{C}_H)=3.
$$
However, it is conjectured that $\alpha_0(\mathcal{S}_H)=5/2$,
which is of special importance in obtaining sharp coefficient estimates for univalent harmonic mappings (see \cite{CS}).
The upper bound for $\alpha_0(\mathcal{S}_H)$ has been improved few times. See \cite[p.10]{CS}, \cite[p.~96]{dur2004} and \cite[Theorem~10]{she}.
However, the conjectured bound remains open. Now the best known upper bound of it was shown in \cite{AAP}.

Let $\mathcal{U}_\alpha$ be the set of all locally univalent analytic functions $h(z)=z+a_2z^2+\cdots$ in $\mathbb{D}$
of order $\leq\alpha$, where $$\alpha=\sup_{h\in\mathcal{U}_\alpha}|a_2(h)|=\frac{1}{2}\sup_{h\in\mathcal{U}_\alpha}|h''(0)|
.$$
The family $\mathcal{U}_\alpha$ is known as the {\it universal linear invariant family} (ULIF) of order $\alpha~(\geq1)$ (see \cite{pom}).
In fact, if $h\in\mathcal{U}_\alpha$, then $K_\varphi(h(z))\in\mathcal{U}_\alpha$ holds for any $\varphi\in{\rm Aut}(\mathbb{D})$.
It is easy to see that a ULIF is not an ALIF, but the set $\{h:\, f=h+\overline{g}\in \mathcal{F}\}$ is a ULIF when $\mathcal{F}$
is an ALIF.
Conversely, we can construct some special ALIFs from ULIFs.
For instance, if
\begin{equation} \label{equ2.1}
\mathcal{F}_\alpha:=\{f: \,f(z)=h(z)+\overline{b_1h(z)},\quad h\in\mathcal{U}_\alpha~~\text{and}~~b_1\in\mathbb{D}\},
\end{equation}
then it is easy to see that $\mathcal{F}_\alpha$ is an ALIF with $\alpha(\mathcal{F}_\alpha)=\alpha_0(\mathcal{F}_\alpha)=\alpha$.
Indeed, simple computation shows that
$$K_\varphi(f(z))=K_\varphi(f(z))+\overline{B_1K_\varphi(f(z))}\in\mathcal{F}_\alpha, \quad B_1=b_1\frac{\varphi'(0)h'(\varphi'(0))}{\overline{\varphi'(0)h'(\varphi'(0))}},
$$
for each $f\in\mathcal{F}_\alpha$ and for any $\varphi\in{\rm Aut}(\mathbb{D})$.
This means that $\mathcal{F}_\alpha$ is a LIF.
On the other hand, calculations prove that $\mathcal{F}_\alpha$ is an AIF because
$$A_\varepsilon(f(z))=h(z)+\left(\frac{\overline{b_1}+\varepsilon}{1+b_1\varepsilon}\right)\overline{h(z)}\in\mathcal{F}_\alpha
$$
for each $f\in\mathcal{F}_\alpha$ and for all $\varepsilon\in\mathbb{D}$.
Moreover, one can also check that $\alpha(\mathcal{F}_\alpha)=\alpha_0(\mathcal{F}_\alpha)=\alpha$.

\subsection{Harmonic pre-Schwarzian derivatives}  \label{sec2.3}

Let $f$ be a locally univalent harmonic mapping in a domain $D$.
The {\it harmonic pre-Schwarzian derivative} of $f$ and the {\it harmonic pre-Schwarzian norm} of $f$ in $D$ are defined by
\begin{equation} \label{equ2.2}
P_f=(\log J_f)_z
\quad \text{and}\quad  ||P_{f}||_D=\sup_{z\in D}\lambda_D^{-1}(z)|P_{f}(z)|,
\end{equation}
respectively (see \cite{HM2015}).
These definitions coincide with the corresponding definitions in the analytic case.

The harmonic pre-Schwarzian derivative inherits the same {\it chain rule} as in the analytic case. More precisely,
if $f$ is a sense-preserving harmonic mapping and $\phi$ is a locally univalent analytic function
for which the composition $f\circ \phi$ is well defined, then, because $J_{f\circ\phi}=|\phi'|^2J_f(\phi)$, we have
\begin{equation}\label{equ2.3}
P_{f\circ\phi}=P_f(\phi)\phi'+P_\phi.
\end{equation}

The harmonic pre-Schwarzian derivative is invariant under an affine transformation of harmonic mapping $f$:
\begin{equation*}
P_{A\circ f}=P_f,\quad A(z)=az+b\overline{z}+c,\quad |a|\neq|b|.
\end{equation*}
So does the harmonic pre-Schwarzian norm.

Besides these, several of the recent properties on this topic may be found from
\cite{gra2016,HM2015,LS2017-1,LS2017-2} and from the references therein.
In these articles, the authors focussed mainly on sense-preserving harmonic mappings $f=h+\overline{g}$ in a simply connected domain $D$
so that the first formulation in \eqref{equ2.2} can be rewritten as
\begin{equation}\label{equ2.4}
P_f=P_h-\frac{\overline{\omega}\omega'}{1-|\omega|^2},
\end{equation}
where $P_h=h''/h'$ and $\omega=\omega_f=g'/h'$.

\section{Inequalities related to harmonic pre-Schwarzian derivatives}  \label{sec3}

In this section, we will generalize results about \eqref{equ1.1}, \eqref{equ1.2},
Theorems \ref{thmA} and \ref{thmB} to univalent harmonic mappings.
However, we will give a negative answer to the case of  \eqref{equ1.3} when $F$ is a univalent harmonic mapping in $\Delta$ with $F(\infty)=\infty$.
The following lemma may be considered as an adjustment of the result of Graf \cite[Theorem~1]{gra2016} as can be
seen by a verification of the corresponding proof. So we omit the details.

\begin{lem} \label{lem1}
Let $\mathcal{F}$ be an {\rm ALIF} and $f\in\mathcal{F}$.
Then
\begin{equation}\label{equ3.1}
|(1-|z|^2)P_{f}(z)-2\overline{z}|\leq2\alpha_0(\mathcal{F}),\quad z\in\mathbb{D},
\end{equation}
and
\begin{equation*}
||P_{f}||_\mathbb{D}=\sup_{z\in\mathbb{D}}(1-|z|^2)|P_{f}(z)|\leq2(\alpha_0(\mathcal{F})+1).
\end{equation*}
Both estimates are sharp if $\mathcal{F}=\mathcal{F}_\alpha$ ($\mathcal{K}_H$, $\mathcal{C}_H$),
where $\mathcal{F}_\alpha$ is defined by \eqref{equ2.1}.
\end{lem}

Recall that $\alpha_0(\mathcal{F}_\alpha)=\alpha$. Let
$f_\alpha= k_\alpha+\overline{b_1k_\alpha}$, where $b_1\in\mathbb{D}$ and $k_\alpha$ is defined by
\be \label{equ3.2}k_\alpha(z)=\frac{1}{2\alpha}\left[1-\left(\frac{1-z}{1+z}\right)^\alpha\right],\quad z\in\mathbb{D}.
\ee
It is easy to see that $f_\alpha\in\mathcal{F}_\alpha$ and \eqref{equ3.1} is sharp for $f_\alpha$ at $z=0$.
Moreover, we have $||P_{f_\alpha}||_{\mathbb{D}}=||P_{k_\alpha}||_{\mathbb{D}}=2(\alpha+1).$
The sharpness part of \eqref{equ3.1} in the cases of $\mathcal{K}_H$ and $\mathcal{C}_H$  can be obtained
by setting $z=0$ and by choosing the harmonic half-plane mapping $L$
and the harmonic Koebe mapping $K$, respectively.
Note that
$$P_L(z)=\frac{3}{1-z}-\frac{\overline{z}}{1-|z|^2}
\quad \mbox{and} \quad P_K(z)=\frac{5+3z}{1-z^2}-\frac{\overline{z}}{1-|z|^2}.
$$
Simple computation and analysis show that $||P_{L}||_\mathbb{D}=5$ (see \cite[Theorem~4]{HM2015})
and $||P_{K}||_\mathbb{D}=7$ (see \cite[Theorem~1.1]{LS2017-1}).
Since the harmonic pre-Schwarzian derivative preserves affine invariance and
$\mathcal{S}_H$ is a special ALIF, we can  get  the following result as a corollary to Lemma \ref{lem1}.
This will be used in the sequel.

\begin{cor} \label{cor1}
Let $f$ be a sense-preserving and univalent harmonic mapping in $\mathbb{D}$.
Then
\begin{equation}\label{equ3.3}
|(1-|z|^2)P_{f}(z)-2\overline{z}|\leq2\alpha_0(\mathcal{S}_H),\quad z\in\mathbb{D},
\end{equation}
and
\begin{equation}\label{equ3.4}
||P_{f}||_\mathbb{D}=\sup_{z\in\mathbb{D}}(1-|z|^2)|P_{f}(z)|\leq2(\alpha_0(\mathcal{S}_H)+1).
\end{equation}
\end{cor}

Clearly, \eqref{equ3.4} is a generalization of  \eqref{equ1.1}.
Now, we will extend \eqref{equ1.2} to the case of harmonic mappings.

\begin{theorem} \label{thm1}
Let $f$ be a sense-preserving and univalent harmonic mapping in
the right half-plane $\mathbb{H}=\{z:\,\RE z>0\}$.
Then
\begin{equation}\label{equ3.5}
||P_{f}||_\mathbb{H}=\sup_{z\in\mathbb{H}}2\RE z|P_{f}(z)|\leq2(\alpha_0(\mathcal{S}_H)+1).
\end{equation}
\end{theorem}
\bpf
Fix $z\in\mathbb{H}$. Then there exists a unique $w\in\mathbb{D}$ such that $z=\phi(w)=\frac{1+w}{1-w}$.
The chain rule \eqref{equ2.3} and a basic computation show that
\begin{align*}
2\RE z|P_{f}(z)|&=2\RE(\phi(w))|P_{f}(\phi(w))|\\
&=\left(\frac{1+w}{1-w}+\frac{1+\overline{w}}{1-\overline{w}}\right)\frac{1}{|\phi'(w)|}|P_{f\circ\phi}(w)-P_{\phi}(w)|\\
&=(1-|w|^2)|P_{f\circ\phi}(w)-P_{\phi}(w)|.
\end{align*}
Note that $f\circ\phi$ is a sense-preserving and univalent harmonic mapping in $\mathbb{D}$ and $P_\phi(w)=2/(1-w)$.
It follows from \eqref{equ3.3} that
\begin{align*}
2\RE z|P_{f}(z)|=&(1-|w|^2)|P_{f\circ\phi}(w)-P_{\phi}(w)| \\
\leq &|(1-|w|^2)P_{f\circ\phi}(w)-2\overline{w}|
+|(1-|w|^2)P_{\phi}(w)-2\overline{w}|\\
\leq&2\alpha_0(\mathcal{S}_H)+2,
\end{align*}
which implies \eqref{equ3.5}.
 \epf

However, there is no similar result to \eqref{equ1.3}
when $F$ is a univalent harmonic mapping in $\Delta=\{z:\,|z|>1\}\cup\{\infty\}$ with $F(\infty)=\infty$.
To demonstrate this fact, we consider
\begin{equation*}
F(z)=z-\frac{1}{\overline{z}}+2\log|z|,\quad z\in\Delta.
\end{equation*}
It follows from \cite[Theorem~3.7]{HS} that $F$ is a sense-preserving and univalent harmonic mapping in $\Delta$.
Direct computations reveal that
$$J_F(z)=\frac{|1+z|^2(|z|^2-1)}{|z|^4}
\quad \mbox{and} \quad P_F(z)=(\log J_F(z))_z=\frac{2+z-|z|^2}{z(z+1)(|z|^2-1)}.
$$
Obviously,
$$\sup_{z\in\Delta}(|z|^3-|z|)|P_F(z)|=\infty.
$$

Next we consider the case of simply connected domain.

\begin{theorem} \label{thm2}
Let $f$ be a sense-preserving and univalent harmonic mapping in a simply connected domain $D\subset\mathbb{C}$.
Then
\begin{equation}\label{equ3.6}
\lambda_D^{-1}(z)|P_{f}(z)|\leq2(\alpha_0(\mathcal{S}_H)+2),\quad z\in D.
\end{equation}
\end{theorem}

\bpf
Fix $z\in D$ and choose a conformal mapping $\phi$ of $\mathbb{D}$ onto $D$ with $\phi(0)=z$.
Distortion theorem of Bieberbach for the univalent function $\phi$ gives that
\begin{equation*}
\left|(1-|\zeta |^2)\frac{\phi''(\zeta )}{\phi'(\zeta )}-2\overline{\zeta}\right|\leq4,\quad \zeta\in\mathbb{D},
\end{equation*}
which implies $|P_\phi(0)|=|\phi''(0)/\phi'(0)|\leq4$.
Moreover, by \eqref{equ1.4}, we have $\lambda_D^{-1}(z)=|\phi'(0)|$.
Note that $f\circ\phi$ is a sense-preserving and univalent harmonic mapping in $\mathbb{D}$,
and thus, \eqref{equ3.3} holds for $f\circ\phi$.
In particular,  \eqref{equ3.3} applied to $f\circ\phi$ at the point $z=0$ yields
$$|P_{f\circ\phi}(0)|\leq 2\alpha_0(\mathcal{S}_H).
$$
By \eqref{equ2.3}, we have
\begin{align*}
\lambda_D^{-1}(z)|P_{f}(z)|&=|\phi'(0)P_{f}(\phi(0))|=|P_{f\circ\phi}(0)-P_\phi(0)|\\
&\leq|P_{f\circ\phi}(0)|+|P_\phi(0)|\\
& \leq 2\alpha_0(\mathcal{S}_H)+4
\end{align*}
and the desired inequality \eqref{equ3.6} follows.
\epf

As remarked in the introduction about  Theorem \ref{thmB}, Theorem \ref{thm2}
cannot be extended to  arbitrary domains.
However, using the method of proof of Theorem \ref{thmB} (we omit the detail),
we obtain a similar result (see Theorem \ref{thm3} below) based on the following lemma,
which is a generalization of \cite[Lemma~1]{osg}.

\begin{lem} \label{lem2}
If $D$ is a proper subdomain of $\mathbb{C}$ and if $f$ is a sense-preserving and univalent harmonic mapping in $D$, then
\begin{equation}\label{equ3.7}
|P_{f}(z)|\leq\frac{2\alpha_0(\mathcal{S}_H)}{d(z,\partial D)},\quad z\in D,
\end{equation}
where $d(z,\partial D)$ denotes the Euclidean distance to the boundary.
\end{lem}
\bpf
Fix $z_0\in D$ and $d_0=d(z_0,\partial D)$.
Then  $F(z)=f(d_0z+z_0)$ is a sense-preserving and univalent  harmonic mapping in $\mathbb{D}$.
Using \eqref{equ2.3} and \eqref{equ3.3}, we get
$$|P_{f}(z_0)|=\frac{|P_{F}(0)|}{d_0}\leq\frac{2\alpha_0(\mathcal{S}_H)}{d_0}
=\frac{2\alpha_0(\mathcal{S}_H)}{d(z,\partial D)}
$$
and the proof is complete.\epf

\begin{theorem} \label{thm3}
Let $D\subset\mathbb{C}$ have at least two boundary points. Then there exists a constant $a$ such that
$\lambda_{D}^{-1}(z)|P_f(z)|\leq a$ in $D$ for all sense-preserving and univalent harmonic mappings in $D$
if and only if there exists a positive constant $c$ such that
\begin{equation*}
\lambda_D(z)d(z,\partial D)\geq c,\quad z\in D.
\end{equation*}
\end{theorem}

\begin{rems}
The inequalities \eqref{equ3.3}-\eqref{equ3.7} are concerned with the conjecture on $\alpha_0(\mathcal{S}_H)$. For example,
if $\alpha_0(\mathcal{S}_H)=5/2$, then these inequalities are sharp. In what follows, $K$ denotes the  harmonic Koebe function.
\begin{enumerate}
\item[(1)] The inequality \eqref{equ3.3} is sharp for the function $K$ (Set $z=0$ in \eqref{equ3.3}).

\item[(2)] The sharpness of \eqref{equ3.4} can be easily seen from the function $K$.

\item[(3)] The sharpness in \eqref{equ3.5} follows by choosing $f=K\circ\phi$ with $\phi(z)=(1-z)/(1+z)$. In order to verify this, we may use \eqref{equ2.2}
and observe that
$$J_f(z)=\frac{z+\overline{z}}{8|z|^8}
~\mbox{ and }~ P_f(z)=-\frac{3z+4 \overline{z}}{z(z+\overline{z})}.
$$

\item[(4)] To show the  sharpness of \eqref{equ3.6}, we  choose $D=\mathbb{C}\backslash(-\infty,0]$ and consider the function
$$f(z)=(K\circ\phi_2\circ\phi_1)(z)=\frac{1}{24} \left(\frac{2}{z^{3/2}}+\frac{3}{z}-5\right)+
\overline{\frac{1}{24} \left(\frac{2}{z^{3/2}}-\frac{3}{z}+1\right)}, \quad z\in D,
$$
where $\phi_1(z)=\sqrt{z}$ ($\arg1=0$) and $\phi_2(z)=\frac{1-z}{1+z}$.
Note that $\phi_1$ (resp. $\phi_2$) is a conformal mapping of $D$ (resp. $\mathbb{H}$)
onto $\mathbb{H}$ (resp. $\mathbb{D}$).
Since $K$ is univalent in $\mathbb{D}$, $f$ is univalent in $D$.
By \eqref{equ2.2}, straightforward computations assert that
$$J_f(z)=\frac{\sqrt{z}+\overline{\sqrt{z}}}{32|z|^5}
\quad \mbox{and} \quad P_f(z)=-\frac{5 \left| z\right| +4 z}{2 z (\left| z\right| +z)}, \quad z\in D.
$$
Note that $\lambda_{D}^{-1}(z)=4z$ and $\lambda_{D}^{-1}(x)|P_f(x)|=9$ for all $x>0$.

\item[(5)] The choice $D=\mathbb{D}$ and $f=K$ show that
\eqref{equ3.7} is also sharp.
\end{enumerate}
\end{rems}

\section{Inequalities between analytic pre-Schwarzian  and harmonic pre-Schwarzian} \label{sec4}

In this section, we try to find some connections between univalent analytic functions and univalent harmonic mappings.
For this, we first investigate certain relationships between the class of analytic pre-Schwarzian derivatives and the harmonic
pre-Schwarzian derivative of a given sense-preserving harmonic mapping.

\begin{theorem} \label{thm4}
Let $f=h+\overline{g}$ be a sense-preserving harmonic mapping in a simply connected domain $D\subset\mathbb{C}$
with the dilatation $\omega_f$.
Then for each $\varepsilon\in\overline{\mathbb{D}}$ we have 
\be \label{equ4.1}
\lambda_D^{-1}(z) |P_{h+\varepsilon g}(z)-P_f(z)|\leq\sup_{z\in D}|\omega_f(z)|,\quad z\in D.
\ee
Moreover, either $||P_{h+\varepsilon g}||_D=||P_f||_D=\infty$ or both $||P_{h+\varepsilon g}||_D$ and $||P_f||_D$
are finite.  
If $||P_f||_D<\infty$, then the inequality
\be \label{equ4.2}
\big| ||P_{h+\varepsilon g}||_D-||P_f||_D \big|\leq\sup_{z\in D}|\omega_f(z)|
\ee
holds.
Furthermore, for any given $k\in[0,1]$, there exists an $\varepsilon_k\in\overline{\mathbb{D}}$ and
a sense-preserving harmonic mapping $F_k=H_k+\overline{G_k}$ in $D$
with  $\sup_{z\in D}|\omega_{F_k}(z)|=k$ such that \eqref{equ4.1} is sharp,
and \eqref{equ4.2} is sharp when $D$ is a disk.
\end{theorem}

\bpf
Let $f=h+\overline{g}$ be sense-preserving in $D$.
Then its dilatation $\omega:=\omega_f:\, D\rightarrow\mathbb{D}$ is analytic
and $k=\sup_{z\in D}|\omega(z)|$  exists, where $k\in [0,1]$.
The proof is trivial for $k=0$ and thus, we assume that $k\in(0,1]$.
For each $\varepsilon\in\overline{\mathbb{D}}$, we observe that
$$|h'(z)+\varepsilon g'(z)|\geq|h'(z)|-|g'(z)|>0,\quad z\in D,
$$
so that $h+\varepsilon g$ is locally univalent in $D$.
Fix $\varepsilon\in\overline{\mathbb{D}}$.
From \eqref{equ2.4}, direct computation shows that
$$P_{h+\varepsilon g}=\frac{h''+\varepsilon g''}{h'+\varepsilon g'}
=P_h+\frac{\varepsilon \omega'}{1+\varepsilon \omega}
$$
and thus,
$$P_{h+\varepsilon g}-P_f=\frac{\varepsilon \omega'}{1+\varepsilon \omega}
+\frac{\overline{\omega}\omega'}{1-|\omega|^2}
=\frac{\varepsilon+\overline{\omega}}{1+\varepsilon\omega}\cdot\frac{\omega'}{1-|w|^2}.
$$

We first consider the special case $D=\mathbb{D}$.
Clearly, $\sup_{z\in\mathbb{D}}\left|\frac{\overline{\varepsilon}
+\omega(z)}{1+\varepsilon\omega(z)}\right|\leq1$.  
On the other hand, applying Schwarz-Pick lemma to  the function $\omega/k: \,\mathbb{D}\rightarrow\mathbb{D}$,
we infer that
$$\frac{|\omega'(z)|(1-|z|^2)}{1-|\omega(z)|^2}
\leq k\frac{1-|\frac{\omega(z)}{k}|^2}{1-|\omega(z)|^2}\leq k,\quad z\in\mathbb{D}.
$$
This means that \eqref{equ4.1} holds for $D=\mathbb{D}$.
Using the triangle inequality, it follows that
$$(1-|z|^2)\big||P_{h+\varepsilon g}(z)|-|P_f(z)|\big|
\leq(1-|z|^2)\left|P_{h+\varepsilon g}(z)-P_f(z)\right|\leq k,\quad z\in\mathbb{D}.
$$
By \eqref{equ2.2},  it is easy to see that either $||P_{h+\varepsilon g}||_\mathbb{D}=||P_f||_\mathbb{D}=\infty$
or both $||P_{h+\varepsilon g}||_\mathbb{D}$ and $||P_f||_\mathbb{D}$ are finite. 
Moreover, if $||P_f||_\mathbb{D}<\infty$, then \eqref{equ4.2} can be deduced from the above inequality when $D=\mathbb{D}$.

Next we will discuss the sharpness part.
For any given $k\in(0,1]$, consider the harmonic function $f_k$ defined on  $\mathbb{D}$ by
\begin{equation} \label{equ4.3}
f_{k}(z)=h_{k}(z)+\overline{g_{k}(z)}
=\int_0^z\frac{(1+kt)^a}{(1-kt)^{a+1}}dt+\overline{k\int_0^z t\frac{(1+kt)^a}{(1-kt)^{a+1}}dt},
\end{equation}
where $a\geq0$.
Obviously,  $\omega_{f_{k}}(z)=kz$ 
and thus, $f_{k}$ is sense-preserving in $\mathbb{D}$.
By computations, we have
$$P_{h_{k}+g_{k}}(z)=\frac{2k(a+1)}{1-k^2z^2}\quad \text{and}
\quad P_{f_{k}}(z)=k\frac{2a+1+kz}{1-k^2z^2}-\frac{k^2\overline{z}}{1-k^2|z|^2}
$$
so that
$$\left[(1-|z|^2)\big|P_{h_{k}+g_{k}}(z)-P_{f_{k}}(z)|\right]_{z=0}=k,
$$
which shows that \eqref{equ4.1} is sharp when $D=\mathbb{D}$.
On the other hand, since
$$(1-|z|^2)|P_{h_{k}+g_{k}}(z)|\leq(1-|z|^2)\frac{2k(a+1)}{1-k^2|z|^2}\leq 2k(a+1),\quad z\in\mathbb{D},
$$ and
$$\left[(1-|z|^2)|P_{h_{k}+g_{k}}(z)|\right]_{z=0}=2(a+1)k,
$$
we see that $||P_{h_{k}+g_{k}}||_\mathbb{D}=2(a+1)k$.
Similarly, we obtain $||P_{h_{k}-g_{k}}||_\mathbb{D}=2ak$.
It follows from \eqref{equ4.2} that $$||P_{h_{k}+g_{k}}||_\mathbb{D}-k\leq||P_{f_{k}}||_\mathbb{D}
\leq||P_{h_{k}-g_{k}}||_\mathbb{D}+k,
$$
which means that $$||P_{h_{k}+g_{k}}||_\mathbb{D}-k=||P_{f_{k}}||_\mathbb{D}
=(2a+1)k=||P_{h_{k}-g_{k}}||_\mathbb{D}+k.
$$
This certifies the sharpness of \eqref{equ4.2} when $D=\mathbb{D}$.

Now we need to consider the general case.
Fix $z\in D$ and consider a conformal mapping $\phi$ of $\mathbb{D}$ onto $D$ with $\phi(0)=z$.
For simplicity, let $f_\varepsilon=h+\varepsilon g$.
Applying \eqref{equ1.4} and \eqref{equ2.3}, we have
\begin{align*}
\lambda_D^{-1}(z)|P_{f_\varepsilon}(z)-P_f(z)|
=&\lambda_D^{-1}(z)|P_{(f_{\varepsilon}\circ\phi)\circ\psi}(z)-P_{(f\circ\phi)\circ\psi}(z)|\\
=&|\psi'(z)|\lambda_D^{-1}(z)|P_{f_{\varepsilon}\circ\phi}(\psi(z))-P_{f\circ\phi}(\psi(z))|\\
=&\lambda_{\mathbb{D}}^{-1}(\psi(z))|P_{f_{\varepsilon}\circ\phi}(\psi(z))-P_{f\circ\phi}(\psi(z))|,
\end{align*}
where $\psi=\phi^{-1}:\, D\rightarrow\mathbb{D}$ is the inverse function of $\phi$.
Thus, \eqref{equ4.1} follows easily because $f\circ\phi$ is a sense-preserving harmonic mapping in $\mathbb{D}$.
To show the sharpness of \eqref{equ4.1}, it suffices to consider the function $F_k=f_{k}\circ\psi$,
where $f_{k}$ is defined by \eqref{equ4.3}.

If $D$ is a disk, then, without loss of generality, we may assume $D=\{z:\, |z-z_0|<r\}~(r>0)$.
The sharpness of \eqref{equ4.2} can be seen from
$$F_k(z)=f_{k}\circ\varphi(z)=H_k(z)+\overline{G_k(z)},\quad z\in D,
$$
where  $f_{k}$ is defined by \eqref{equ4.3} and $\varphi(z)=(z-z_0)/r$.
In fact, it follows from \eqref{equ1.4} and \eqref{equ2.3} that
\begin{align*}
||P_{F_k}||_D=&\sup_{z\in D}\lambda_D^{-1}(z)|P_{F_k}(z)|=\sup_{z\in D}|\varphi'(z)|\lambda_D^{-1}(z)|P_{f_{k}}(\varphi(z))|\\
=&\sup_{z\in D}\lambda_\mathbb{D}^{-1}(\varphi(z))|P_{f_{k}}(\varphi(z))|
=\sup_{z\in \mathbb{D}}\lambda_\mathbb{D}^{-1}(z)|P_{f_{k}}(z)|\\
=&||P_{f_{k}}||_\mathbb{D}=(2a+1)k.
\end{align*}
By a similar analysis, we get that
$$||P_{H_k+G_k}||_D-k=||P_{F_{k}}||_D=(2a+1)k=||P_{H_k-G_k}||_D+k,
$$
which completes the proof.
\epf

Applying the triangle inequality twice to \eqref{equ4.1} and \eqref{equ4.2},
we acquire the following sharp results.
Furthermore, the sharpness follow from the same lines of the proof of Theorem \ref{thm4}.

\begin{cor} \label{cor2}
Let $f=h+\overline{g}$ be a sense-preserving harmonic mapping in a simply connected domain $D\subset\mathbb{C}$
with the dilatation $\omega_f$.
Then for each pair $\varepsilon_1,~\varepsilon_2\in\overline{D}$,  we have 
\be \label{equ4.4}
\lambda_D^{-1}(z)\big| P_{h+\varepsilon_1 g}(z)-P_{h+\varepsilon_2 g}(z)\big|
\leq2\sup_{z\in D}|\omega_f(z)|, \quad z\in D.
\ee
Moreover, either $||P_{h+\varepsilon_1 g}||_D=||P_{h+\varepsilon_2 g}||_D=\infty$
or both $||P_{h+\varepsilon_1 g}||_D$ and $||P_{h+\varepsilon_2 g}||_D$ are finite.
If $||P_{h+\varepsilon g}||_D<\infty$ for some $\varepsilon\in\overline{\mathbb{D}}$, then
\be \label{equ4.5}
\big| ||P_{h+\varepsilon_1 g}||_D-||P_{h+\varepsilon_2 g}||_D \big|\leq2\sup_{z\in\mathbb{D}}|\omega_f(z)|
\ee
holds for any $\varepsilon_1,~\varepsilon_2\in\overline{\mathbb{D}}$.
Furthermore, for any given $k\in[0,1]$, there exist $\varepsilon_1(k),~\varepsilon_2(k)\in\overline{\mathbb{D}}$ and
a sense-preserving harmonic mapping $F_k=H_k+\overline{G_k}$ in $D$
with $\sup_{z\in D}|\omega_{F_k}(z)|=k$ such that \eqref{equ4.4} is sharp,
and \eqref{equ4.5} is sharp when $D$ is a disk.
\end{cor}

\begin{cor} \label{cor3}
For any sense-preserving harmonic mapping $f=h+\overline{g}$ in a simply connected domain $D\subset\mathbb{C}$,
we have
\be \label{equ4.6}
\lambda_D^{-1}(z) |P_{h+\varepsilon g}(z)-P_f(z)|\leq1 \quad \forall~\varepsilon\in\overline{\mathbb{D}}
\ee
in $D$ and
\be \label{equ4.7}
\max\left\{0,~\max_{\varepsilon\in \overline{\mathbb{D}}} ||P_{h+\varepsilon g}||_D-1\right\}
\leq||P_f||_D\leq \min_{\varepsilon\in \overline{\mathbb{D}}} ||P_{h+\varepsilon g}||_D+1.
\ee
The constant 1 is sharp in the two estimates.
\end{cor}

In particular, if $f=h+\overline{g}$ in Corollary \ref{cor3} is further restricted to be univalent, then,
for any given $\varepsilon\in\overline{\mathbb{D}}$,  the distance between $||P_f||_D$ and $||P_{h+\varepsilon g}||_D$ is at most 1.
Compared to the corresponding results in Sections \ref{sec1} and \ref{sec3},
if the  conjecture of Clunie and Shell-Small \cite{CS} on $\alpha_0(\mathcal{S}_H)$ were true,
then the distance between the sharp constant 6 (resp. 8) in \eqref{equ1.1} and \eqref{equ1.2} (resp. \eqref{equ1.5})
 and the sharp constant $2(\alpha_0(\mathcal{S}_H)+1)$ (resp. $2(\alpha_0(\mathcal{S}_H)+2)$)
 in \eqref{equ3.4} and \eqref{equ3.5} (resp. \eqref{equ3.6}) is also 1.
This raises the following conjecture.

\begin{conjecture} \label{conj1}
Let  $f=h+\overline{g}$ be a sense-preserving and univalent harmonic mapping in a simply connected domain $D\subset\mathbb{C}$.
Then there exists   a constant $\varepsilon\in\overline{\mathbb{D}}$ such that
$h+\varepsilon g$ is univalent in $D$.
\end{conjecture}

It is easy to see that to solve the above conjecture, it suffices to consider the case $D=\mathbb{D}$.
Moreover, this conjecture is weaker than the following conjecture proposed in \cite{PS}.

\begin{conj} \label{conjA}
For every function $f=h+\overline{g}\in\mathcal{S}_H^0$, there exists  a constant $\theta\in\mathbb{R}$ such that
$h+e^{i\theta} g\in\mathcal{S}$, where $\mathcal{S}$ is the class of analytic functions in $\mathcal{S}_H$.
\end{conj}

If Conjecture \ref{conj1} were true, then, combining \eqref{equ4.6} with Bieberbach's distortion theorem for univalent analytic functions in $\mathbb{D}$,
we can obtain the sharp inequality \eqref{equ3.3}.
Moreover, combining \eqref{equ4.7} (resp. \eqref{equ4.6}) with \eqref{equ1.1} and \eqref{equ1.2} (resp. \eqref{equ1.5}),
we can obtain sharp inequalities \eqref{equ3.4} and \eqref{equ3.5} (resp. \eqref{equ3.6}), respectively.
In fact, Conjecture \ref{conj1} implies $\alpha_0(\mathcal{S}_H)=5/2$. 
To clarify this, let us assume for the moment that Conjecture \ref{conj1} is true. Then, for  $f=h+\overline{g}\in\mathcal{S}_H^0$,
there exists an $\varepsilon_0\in\overline{\mathbb{D}}$ such that
$h+\varepsilon_0 g$ is univalent in $\mathbb{D}$.
Let $D=\mathbb{D}$ in \eqref{equ4.6}. Then we have that $|P_{h+\varepsilon_0g}(0)-P_f(0)|\leq1$.
Using the method of proof of Theorem \ref{thm2}, we can obtain that $|P_{h+\varepsilon_0g}(0)|\leq4$.
Note that $|P_f(0)|=|P_h(0)|=2|a_2(f)|$, because of $\omega_f(0)=0$.
From this observation, it is easy to see that $|a_2(f)|\leq5/2$ holds.
It is well known that the harmonic Koebe function $K$ belongs to $\mathcal{S}_H^0$ with $a_2(K)=5/2$,
which implies that $\alpha_0(\mathcal{S}_H)=5/2$.

Note that Conjecture \ref{conjA} not only implies Conjecture \ref{conj1}, but also provides an affirmative answer to the open question
about the estimates of the coefficients for harmonic mappings in $\mathcal{S}_H^0$ (see \cite{CS}).

\section{Applications} \label{sec5}

In this section, we improve the corresponding results of \cite[section~3]{sch}
where the author used the ALIF of harmonic mappings to obtain the radius of majorization of
the Jacobian of harmonic mappings.
A function $f$ is said to be majorized by $F$ in a certain region if $|f(z)|\leq|F(z)|$ holds there.
Next let's recall the notion of subordination.
Let $f$ and $F$ be two harmonic mappings in  $\mathbb{D}$.
We say that $f$ is subordinate to $F$, denoted by $f\prec F$, if $f(z)=F(\psi(z))$,
where $\psi$ is analytic with $\psi(0)=0$ and $|\psi(z)|<1$ in $\mathbb{D}$.
For  the importance, background, development and results concerning these two topics,
the reader may refer to the paper \cite{sch}  and  the references therein.

Below, we denote $n(x)=1+x-\sqrt{x^2+2x}$, $x\geq0$.
One of the subordination-majorization results for ULIF of analytic functions is the following.

\begin{thm}  \label{thmC} 
Let $\mathcal{U}_\alpha$ be a {\rm ULIF} with $1\leq\alpha(\mathcal{U}_\alpha)<\infty$.
If  $f\prec F$ and $F\in\mathcal{U}_\alpha$, then $|f'(z)|\leq |F'(z)|$ for
$|z|\leq n(\alpha)$, $\alpha=\alpha(\mathcal{U}_\alpha)$ and the result is best possible.
\end{thm}

The case $\alpha(\mathcal{U}_\alpha)\geq1.65$ in Theorem \ref{thmC} was  established first by Campbell (see \cite{cam}).
Also, he conjectured that Theorem \ref{thmC} held for $1\leq\alpha(\mathcal{U}_\alpha)<1.65$,
which was later proved affirmatively  by Barnard and Pearce (see \cite{BP}).
By the way, the case $\alpha(\mathcal{U}_\alpha)=1$ was dealt in \cite{BK}.
It is natural to ask for the analogous result to ALIF of harmonic mappings.
Schaubroeck \cite{sch} has obtained a partial answer.
We will improve his results based on the following theorem,
 which is exactly the same as \cite[Theorem~1]{gra2012}, but we present a much simpler proof of it.

\begin{theorem} \label{thm5}
Let $\mathcal{F}$ be an {\rm ALIF} and $f\in\mathcal{F}$. Then
\be \label{equ5.1}
\frac{(1-r)^{2\alpha-2}}{(1+r)^{2\alpha+2}}\leq\frac{J_f(z)}{1-|b_1|^2}
\leq\frac{(1+r)^{2\alpha-2}}{(1-r)^{2\alpha+2}},\quad |z|=r<1,
\ee
where $\alpha=\alpha_0(\mathcal{F})$ and $b_1=\overline{f_{\overline{z}}(0)}$.
Equalities occur if $\mathcal{F}=\mathcal{F}_\alpha$
for the functions $f(z)=k_\alpha(z)+\overline{b_1k_\alpha(z)}$,
where $\mathcal{F}_\alpha$ and $k_\alpha$ are defined by \eqref{equ2.1} and \eqref{equ3.2}, respectively.
In addition, equalities hold if  $\mathcal{F}=\mathcal{K}_H$ (resp. $\mathcal{C}_H$)
for the functions
$$f(z)=L(z)+\overline{b_1L(z)}\quad \left(\text{resp.} ~K(z)+\overline{b_1K(z)}\right),
$$
where $L$ and $K$ are the harmonic half-plane mapping and the harmonic Koebe mapping, respectively.
\end{theorem}
\bpf
Let $z=re^{i\theta}\in\mathbb{D}$.
It follows from \eqref{equ3.1} that $|z(1-|z|^2)P_f(z)-2|z|^2|\leq2\alpha|z|$ and thus
\begin{equation} \label{equ5.2}
\frac{2r^2-2\alpha r}{1-r^2}\leq \RE zP_f(z)\leq\frac{2\alpha r+2r^2}{1-r^2},\quad |z|=r<1.
\end{equation}
Note that
$$\frac{1}{2}r\frac{\partial}{\partial r}\log J_f(z)
=\frac{r}{2}(P_f(z)e^{i\theta}+\overline{P_f(z)}e^{-i\theta})=\RE zP_f(z).
$$
If we substitute the above equality into \eqref{equ5.2} and integrate the resulting inequalities, we obtain the desired conclusion.
\epf

The following result is a generalization of Theorem \ref{thmC}.
The proof of this result follows if we adopt the method of the proof of \cite[Theorem~3.4]{sch} carefully.
For the sake of completeness, we include the details here.

\begin{theorem} \label{thm6}
Let $\mathcal{F}$ be an {\rm ALIF} with $1\leq\alpha_0(\mathcal{F})<\infty$.
If $f\prec F$ and $F\in\mathcal{F}$, then $J_f(z)\leq J_F(z)$ for
$|z|\leq n(\alpha_0(\mathcal{F}))$ and the result is best possible if  $\mathcal{F}=\mathcal{F}_\alpha$,
where $\mathcal{F}_\alpha$ is defined by \eqref{equ2.1}.
\end{theorem}
\bpf Fix $F=H+\overline{G}\in\mathcal{F}$  and consider
$$S(F)(z)=A_{\varepsilon_0}\circ K_\varphi(F(z)), \quad \varphi(z)=\frac{z_0+z}{1+\overline{z_0}z},~
\varepsilon_0=-\frac{\overline{G'(\varphi(0))}}{H'(\varphi(0))},
$$
where $|z_0|<1$, $K_\varphi$ and $A_{\varepsilon_0}$ are defined as in Section \ref{sec2.2}.
It is easy to see that $S(F)\in\mathcal{F}^0$. Denote $\alpha_0(\mathcal{F})$ by $\alpha$.
It follows from \eqref{equ5.1} that
\begin{equation} \label{equ5.3}
J_{S(F)}(z)=\frac{|\varphi'(z)|^2}{(1-|z_0|^2)^2}\cdot \frac{J_F(\varphi(z))}{J_F(z_0)}
\leq\frac{(1+|z|)^{2\alpha-2}}{(1-|z|)^{2\alpha+2}}.
\end{equation}
If we substitute $z=(x-z_0)/(1-\overline{z_0}x)$ into \eqref{equ5.3}, then a direct calculation shows that
\begin{equation} \label{equ5.4}
\frac{J_F(x)}{J_F(z_0)}\leq\left(\frac{1-|z_0|^2}{1-|x|^2}\right)^2
\left(\frac{1+|(x-z_0)/(1-\overline{z_0}x)|}{1-|(x-z_0)/(1-\overline{z_0}x)|}\right)^{2\alpha}.
\end{equation}
If $f\prec F$, then $f(z)=F(\psi(z))$ and thus, $J_f(z)=|\psi'(z)|^2J_F(\psi(z))$.
Applying \eqref{equ5.4} at $x=\psi(z_0)$, we obtain that
$$\sqrt{\frac{J_f(z_0)}{J_F(z_0)}}\leq|\psi'(z_0)|\frac{1-|z_0|^2}{1-|\psi(z_0)|^2}
\left(\frac{|1-\overline{z_0}\psi(z_0)|+|\psi(z_0)-z_0|}{|1-\overline{z_0}\psi(z_0)|-|\psi(z_0)-z_0|}\right)^{\alpha}.
$$
The quantity on the right side above is not greater than $1$ for $|z|\leq n(\alpha)$ (see  \cite{BP} and \cite{cam}).

To show that the result is best possible for  $\mathcal{F}=\mathcal{F}_\alpha$,
it suffices to consider the following functions
$$F(z)=k_\alpha(z)+\overline{b_1k_\alpha(z)}\quad \text{and}\quad f_a(z)=F(\psi(z)):= h_a(z)+\overline{b_1h_a(z)},
$$
where $b_1\in\mathbb{D}$, $k_\alpha$ is defined by \eqref{equ2.2} and $\psi(z)=z(a+z)/(1+az)$, $0\leq a\leq 1$.
Clearly, $f_a(z)\prec F(z)$ for any $a\in[0,1]$ and $F\in\mathcal{F}_\alpha$.
Note that
$$h_a(z)=k_\alpha(\psi(z)),~ J_F(z)=(1-|b_1|^2)|k'_\alpha(z)|^2~ \mbox{ and }~ J_{f_a}(z)=(1-|b_1|^2)|h'_a(z)|^2.
$$
It follows from the proof of \cite[p.302-303]{cam} that the result in Theorem \ref{thmC} is best possible.
Therefore, the result in Theorem \ref{thm6} is also best possible if  $\mathcal{F}=\mathcal{F}_\alpha$.
\epf

\begin{cor} \label{cor4}
If $f\prec F$ and $F\in \mathcal{K}_H$ (resp. $\mathcal{C}_H$),
then $J_f(z)\leq J_F(z)$ for $|z|\leq n(3/2)=(5-\sqrt{21})/2\approx 0.208712$
(resp. $n(5/2)=(7-3\sqrt{5})/2\approx0.145898$) and the result is best possible.
\end{cor}
\bpf
Note that $\alpha_0(\mathcal{K}_H)=3/2$,  $\alpha_0(\mathcal{C}_H)=5/2$ and thus, it suffices to prove the sharpness for both cases.
We first show that the result for $\mathcal{K}_H$ cannot be improved. For this, we let $F(z)=-L(-z)$ and $f_a(z)=F(\psi(z))$,
where $L$ is the harmonic half-plane mapping and $\psi(z)=z(a+z)/(1+az)$, $0\leq a\leq 1$.
Obviously, $f_a\prec F$  for each $a\in[0,1]$ and $F\in\mathcal{K}_H$.
Direct computations give for $0<r<1$,
$$J_{f_a}(r)=(1-r^2)\frac{(a+2r+ar^2)^2}{(1+2ar+r^2)^5}
\quad \mbox{and}\quad \left[\frac{\partial}{\partial a}J_{f_a}(r)\right]_{a=1}=2\frac{1-r}{(1+r)^7}(1-5r+r^2)
$$
which is negative in the interval $(0,1)$ if $n(3/2)<r<1$. This means that $J_{f_a}(r)$ is a decreasing function of $a$ for $r$ in this interval.
Therefore, for any given sufficient small $\varepsilon>0$, there exists a constant $a_\varepsilon$
which is sufficiently close to $1$ such that
$$J_{f_{a_\varepsilon}}(r)>J_{f_1}(r)=\frac{1-r}{(1+r)^5}=J_F(r),\quad n(3/2)+\varepsilon<r<1.
$$
Thus, $J_{f_{a_\varepsilon}}$ is not majorized by $J_F$ in $n(3/2)+\varepsilon<|z|<1$
which shows that the number $n(3/2)$ is best possible.

Let us next indicate the sharpness of the result for the family $\mathcal{C}_H$. Indeed,
similar analysis for the family $\mathcal{C}_H$ may be considered by setting $F(z)=-K(-z)$ and
$f_a(z)=F(\psi(z))$, where $K$ is the  harmonic Koebe function, $\psi(z)=z(a+z)/(1+az)$ with $0\leq a\leq 1$, and
for $0<r<1$,
$$J_{f_a}(r)=(1-r^2)^3\frac{(a+2r+ar^2)^2}{(1+2ar+r^2)^7}
\quad \mbox{and}\quad \left[\frac{\partial}{\partial a}J_{f_a}(r)\right]_{a=1}=2\frac{(1-r)^3}{(1+r)^9}(1-7r+r^2).
$$
This completes the proof.
\epf
\subsection*{Acknowledgments}
The work was completed during the visit of the first author to the Indian Statistical Institute,
Chennai Centre and this author thanks the institute for the support and the hospitality.
The research of the first author was supported by the NSFs of China~(No. 11571049),
the Construct Program of the Key Discipline in Hunan Province and
the Science and Technology Plan Project of Hunan Province (No. 2016TP1020).
The second author is on leave from IIT Madras, Chennai.


\end{document}